\documentclass{elsart}               

\usepackage{graphicx}          
\usepackage{natbib}
\usepackage{amsmath}
\usepackage{amssymb}
\usepackage{colortbl}

\newcommand{\clr}{\color{black} }
\newcommand{\clb}{\color{black} }

\begin{document}

\begin{frontmatter}

\title{A note on relationship between some bounding inequalities in stability analysis of time-delay systems}

\author[gye]{\'Eva Gyurkovics }\ead{gye@math.bme.hu}, 
\author[tt]{Tibor Tak\'acs }\ead{takacs.tibor@uni-corvinus.hu}

\address[gye]{Mathematical Institute, Budapest University of Technology and Economics, M\H{u}egyetem rkp. 3, Budapest, Hungary}
\address[tt]{Corvinus University of Budapest,
8 F\H{o}v\'am t\'er, H-1093, Budapest, Hungary}


\begin{abstract}                          
In this paper, an extension of the generalized free matrix based inequality is introduced in a unified form suitable for the estimation of integrals and sums of quadratic functions. The equivalences of several known variants are shown, including the free matrix based inequalities and its simplified form. Secondly, the relationship between the (simplified) free matrix based estimation and the combination of the Bessel-based inequality with different bounding inequalities affine in the length of the intervals are investigated.
\end{abstract}

\end{frontmatter}

\section{Introduction}
The stability of time-delay systems is often analyzed by means of appropriate Lyapunov-Krasovskii functionals, with the help of which
tractable stability criteria can be derived (see e.g. \cite{chen16a} - \cite{zha17} and the references therein).
In these investigations, the lower estimation for integrals (in case of continuous time systems) or sums (in case of discrete time systems) of positive quadratic terms plays a crucial role.  A possible tool is the Jensen's inequality, but in the past years several results were published to reduce its conservatism. The so called free-matrix-based (FMB) and generalized free-matrix-based (GFMB) inequalities has been proposed recently for this {\clr purpose} (see {\clr \cite{zeng15a}}, \cite{chen16a}, {\clr \cite{zha16a},} \cite{lee17a}, \cite{lee17b}). Another line of the improvement of Jensen's inequality is represented by the Wirtinger inequality, the Bessel-Legendre and Bessel-Chebyshev inequalities (see e.g. \cite{S-GAut13}, \cite{seu14}, {\clr \cite{hien15}} \cite{gyet16}, \cite{gyet17}). The relationship between the different approaches may be of interest. It has been proven in paper \cite{gye15} that the Wirtinger inequality of \cite{S-GAut13} and the FMB inequality of \cite{zeng15a} are equivalent with superiority of the  Wirtinger inequality. This result has been generalized lately by \cite{chen16b}. However, the Jensen's and Wirtinger'{\clr s} inequalities and their generalizations depend in a non convex way {\clb on} the length of the intervals, thus they can only be applied in combination with a convexifying inequality as e.g. reciprocally convex combination lemma of {\clr \cite{park11} when} time-varying-delay systems are considered. Therefore, such combinations have to be taken into account to obtain a real comparison.
\emph{The main purpose of the present paper} is twofold. First, an extension of the generalized free matrix based inequality is to be formulated in a unified form, 
and the relationship {\clr between its  variants} is to be investigated. 
Secondly, the estimations obtained by FMB 
 and by the combination of the Bessel-based inequality with different bounding inequalities affine in the length of the intervals are to be compared.
Throughout this paper, 
 {\clr  $\mathbf{S}_{n}^{+}$ is } the set of positive definite symmetric matrices of size $n\times n,$ and
 $\mathop{\textrm{He}}(A)=A+A^T ,$ where $A^T$ is the transpose of $A.$
\section{Generalized free-matrix-based approach and its variants}
In this section, a generalized free-matrix-based estimation will be formulated by a lemma that will serve as a base of comparisons of different
related approaches that can be applied for time-delay systems  both in continuous and discrete time cases. In this formulation, we shall keep in mind that, in the application for the stability analysis of systems with time-varying delays, one splits an interval to two subsequent subintervals.
Let $\mathbb{E}_i$ $(i=0,1,2)$ be the Euclidean space of functions $\varphi:D_i \subset \mathbf{R} \rightarrow \mathbf{R}$ with the scalar product
$\langle .,. \rangle _i$ containing the elements $\Pi_{0i}(t)\equiv 1,$ $t \in D_i,$ respectively, and possessing with the following two properties:

\vspace{-0.18cm}
\begin{description}
  \item[(P1)] If $\varphi , \psi \in \mathbb{E}_i ,$ then $\varphi  \psi \in \mathbb{E}_i $ and $\langle \varphi , \psi \rangle_{{\clr i}} = \langle  \Pi_{0i}, \varphi  \psi \rangle_{{\clr i}} ;$
  \item[(P2)] If for $\varphi \in \mathbb{E}_i $ $\varphi(t) \geq 0$ for all $t \in D_i,$ then $\langle \Pi_{0i}, \varphi  \rangle_{{\clr i}} \geq 0 .$
\end{description}

\vspace{-0.18cm}
\begin{rem}\label{rem:0}
Typically, $D_i \subset  \mathbf{R} \mbox{ or } \mathbf{Z}, $ $ D_0=[a,b],$ $ D_1=[a,c),$ and $ D_2=[c,b],$ and the scalar products are {\clr defined as} $\langle \varphi , \psi \rangle_{{\clr i}} =\int_{D_i}\varphi(t)\psi(t)dt$ or $\langle \varphi , \psi \rangle_{{\clr i}} =\sum_{t \in D_i}\varphi(t)\psi(t),$ having properties (P1) and (P2).
\end{rem}

Let $\Pi_{0i}, \Pi_{1i},..., \Pi_{\nu i}$ be an orthogonal system in $\mathbb{E}_i$ for some non-negative integer $\nu ,$ and
consider $f \in \mathbb{E}_{{\clr 0}} ^n .$ Suppose that $D_1 \cap D_2 = \emptyset,$ $D_1 \cup D_2 =D_0$ and $\langle \varphi , \psi \rangle_0= \langle \varphi _1, \psi _1 \rangle_1 +
\langle \varphi _2, \psi _2 \rangle_2 ,$
where $\varphi_i$ and $\psi _i$ are the restrictions of $\varphi$ and $\psi$ respectively to $D_i .$
Set {\clr $M_1 =(\nu +1)n,$  $M_2=2M_1 $} and consider $\underline{w} \in \mathbf{R}^{M_{{\clr 2}}},$ with
\begin{align*}
\underline{w} =
                  \begin{bmatrix}
                    \underline{w}^{1} \\
                    \underline{w}^{2} \\
                  \end{bmatrix}
                      =
                           \begin{bmatrix}
                             \mathop{\textrm{col}} \left\{ \langle f_1,\Pi_{01}\rangle ,..., \langle f_1,\Pi_{\nu 1}\rangle   \right\} \\
                             \mathop{\textrm{col}} \left\{ \langle f_2,\Pi_{02}\rangle ,..., \langle f_2,\Pi_{\nu 2}\rangle   \right\} \\
                           \end{bmatrix}  ,
\end{align*}
where $f_i$ is the restriction of $f$ to $D_i .$ Set furthermore symmetric matrices $\Psi^i $ $(i=1,2)$ as
\begin{equation*}
\Psi^i =
           \begin{bmatrix}
             Z_{00}^{i} & \ldots & Z_{0 \nu}^{i} & N_{0}^{i} \\
             \vdots & \ddots & \vdots & \vdots \\
             Z_{\nu 0}^{i} & \ldots & Z_{\nu 0}^{\nu 0} & N_{\nu}^{i} \\
              {N_{0}^{i}}^T & \ldots & {N_{\nu}^{i}}^T & W \\
           \end{bmatrix}      ,
\end{equation*}
where
$Z_{kl}^{i} \in \mathbf{R}^{M_{{\clr 2}} \times M_{{\clr 2}}} ,$ $ N_{k}^i \in \mathbf{R}^{M_{{\clr 2}} \times n},$ $k,l=0. \ldots,\nu .$
\begin{lem} [GFMB inequality.] \label{lem:1}
If $W \in \mathbf{S}_{n}^{+} ,$  and
\begin{align}\label{c26_1}
  \Psi^i & \geq 0,
\end{align}
then the following generalized free-matrix-based inequality holds true {\clr for constants
$\rho _{k}^{i} \geq \| \Pi_{ki} \|_{i}^{2} $}:
\begin{align}
& \langle f_i , Wf_i \rangle  \geq
   - \chi _{i}^T \left( \sum _{k=0}^{\nu} \rho_{k}^{i} Z_{kk}^{i} \right) \chi _{i} -  \mathop{\emph{\textrm{He}}} \left(   \chi _{i}^{T} N^i \underline{w}, \right) \label{ttt1}
\end{align}
where $\chi _{i} \in \mathbf{R}^{M_{{\clr 2}} }$ ($i=1,2$) is arbitrary ${N}^1 = \left( \widehat{N}^1 ,0 \right)= $ $ \left( N_0^1 ,..., N_{\nu}^{1} ,0,...,0 \right)   \in \mathbf{R}^{M_{{\clr 2}} \times M_{{\clr 2}}}$ and
$N^2 =  \left( 0, \widehat{N}^2 \right) $ $=\left( 0,...,0, N_0^2 ,..., N_{\nu}^{2} \right)  \in \mathbf{R}^{M_{{\clr 2}} \times M_{{\clr 2}}} .$
\end{lem}
\textbf{Proof.} For $i=1,2,$ set $ \xi _{i} = \mathop{\textrm{col}}\left\{ \Pi_{0i} \chi _{i} ,...
, \Pi_{\nu i} \chi _{i}, f_i \right\}.$ Then $\xi _{i} (t) \in \mathbf{R}^{(\nu +1) M_{{\clr 2}} +n}$ and $\xi _{i} (t)^{T} \Psi^i \xi _{i} (t) \geq 0 .$
We obtain that
\begin{align}
0 \leq & \langle  \Pi_{0i}, \xi _{i}^T \Psi^i \xi _{i} \rangle _i = \langle \xi _{i}, \Psi^i \xi _{i}\rangle _i  \nonumber   \\
= &
 \sum_{k=0}^{\nu} \sum_{l=0}^{\nu} \langle \Pi_{ki},\Pi_{li} \rangle _i \ \chi _{i}^T Z_{kl}^{i} \chi _{i} 
 + 2 \sum_{k=0}^{\nu} \chi _{i}^T
N_{k}^{i} \langle f_i,\Pi_{ki} \rangle _i + \langle f_i,W f_i \rangle_i .    \label{ttt2}
\end{align}
Since $\Psi ^{i} \geq 0,$ $Z_{kk}^{i} \geq 0$ is true.
Taking into consideration the orthogonality of $\Pi_{ki}$s and the definition $\langle f_i,\Pi_{ki} \rangle =w_{k}^{i},$ (\ref{ttt1}) can be obtained by rearranging (\ref{ttt2}).

\begin{rem}\label{rem:A0}
Frequently, Lemma \ref{lem:1} is applied for the derivative (in continuous-time case) or the forward difference (in discrete-time case) of a function under the choice of $\left\{ \Pi_{ki} \right\}_{k=0}^{\nu}$ as orthogonal polynomials. Using partial integration, i.e. the Abel-lemma, and expressing the derivatives in a chosen basis of polynomials, one can derive analogous estimations. (See the details how such derivations can be performed e.g. in \cite{seu14}, \cite{gyet16}, \cite{leejfi15}, \cite{chen16e}, etc. for continuous time, and \cite{gyet17}, \cite{hien15}, \cite{zha16b}, etc. for discrete time.)
\end{rem}

\subsection{Independent functions based GFMB inequality}

Recently, \cite{lee17a} and \cite{lee17b} derived a GFMB estimation in Lemma {\clb 3} {\clr and in} Lemma 7 by
considering certain linearly independent system of polynomials instead of orthogonal ones. We shall formulate a unified paraphrase of these lemmas as follows.
\begin{lem} [IFB-GFMB inequality.] \label{lem:2}
Let $\left\{ p_{ki} \right\}_{k=0}^{\nu}$ be a system of linearly independent functions defined on $\mathbb{E}_i$, $\widetilde{w}_{k}^{i} = \langle f_i,p_{ki} \rangle _i$ and
$\gamma_{kl} = \langle p_{ki},p_{li} \rangle_i .$
If $W \in \mathbf{S}_{n}^{+} ,$  and (\ref{c26_1}) holds true,
then the following independent-functions-based generalized free-matrix-based inequality holds true:
\begin{align}
  \langle f_i , W f_i \rangle
\geq &- \sum _{k=0}^{\nu} \mathop{\textrm{He}} \left(  \chi _{i}^T N_{k}^{i} \widetilde{w}_{k}^{i}  \right)
\nonumber  \\
&   %
-  \chi _{i}^T \left(  \sum _{k=0}^{\nu} \gamma _{kk}Z_{kk}^{i} +
\sum _{k=0}^{\nu} \sum _{l=k+1}^{\nu} \mathop{\textrm{He}} \left( \gamma_{kl} Z_{kl}^{i} \right)  \right) \chi _{i} .  \label{ttt3}
\end{align}
where $\chi _{i} \in \mathbf{R}^{M_{{\clr 2}} }$ ($i=1,2$) is arbitrary.
\end{lem}

First we show that, although (\ref{ttt3}) contains more free parameters than (\ref{ttt1}), it does not give a better lower estimation.
\begin{thm}\label{thm:1} Suppose that $\left\{ p_{ki} \right\}_{k=0}^{\nu}$ and $\left\{ \Pi_{ki} \right\}_{k=0}^{\nu}$ span the same subspace of $\mathbb{E}_i$. Let $\rho _{k}^{i} = \| \Pi_{ki} \|_{i}^{2} .$ Then
estimations GFMB and IFB-GFMB are equivalent.
\end{thm}
\textbf{Proof.} Since the orthogonal system of functions is also linearly independent,
it is enough to show that GFMB implies IFB-GFMB. For simplicity, we omit index $i$ in the proof. Let ${\Psi} \geq 0$  be given. We will show that a
$\widetilde{\Psi} \geq 0$ exists such that
\begin{align}
&- \chi ^T \left( \sum _{k=0}^{\nu} \|  \Pi _k \|^{2} \widetilde{Z}_{kk} \right)
\chi  - \mathop{\textrm{He}} \left( \chi ^T \widetilde{N}_k \langle f,\Pi _k \rangle \right) \nonumber  \\
= &- \chi ^T \left( \sum _{k=0}^{\nu} \gamma _{kk}Z_{kk} +
\sum _{k=0}^{\nu} \sum _{l=k+1}^{\nu} \mathop{\textrm{He}} \left( \gamma_{kl} Z_{kl} \right)  \right) \chi 
- \sum _{k=0}^{\nu} \mathop{\textrm{He}} \left( \chi ^T N_{k} \langle f,p_k \rangle  \right) .  \label{ttt4}
\end{align}
In fact, there is an invertible matrix $C \in $ $\mathbf{R}^{(\nu + 1) \times (\nu + 1)}$ such that
$\mathop{\textrm{col}} \left\{ p_0,...,p_{\nu} \right\} $ $= C \mathop{\textrm{col}} \left\{ \Pi_0,...,\Pi_{\nu} \right\},$ i.e.
$p_i = \sum _{j=0}^{\nu} c_{ij} \Pi_j ,$ therefore
\begin{align*}
\gamma_{kl} &= \langle p_k,p_l \rangle
= \sum _{j=0}^{\nu}c_{kj} c_{lj} \| \Pi _j \|^2 .
\end{align*}
Let $\widetilde{\Psi}$ be defined by
$\widetilde{Z}_{jj}=  \sum _{k=0}^{\nu} \sum _{l=0}^{\nu} c_{kj}c_{kl}Z_{kl} ,$
$\widetilde{N}_j =\sum _{k=0}^{\nu} c_{kj}N_k$ and $\widetilde{Z}_{kl}=\widetilde{N}_k W^{-1}\widetilde{N}_{l}^{T},$ if $k \neq l . $
Then we obtain that
\begin{align}
\chi ^T \left( \sum _{k=0}^{\nu} \sum _{l=0}^{\nu} \gamma_{kl} Z_{kl} \right) \chi &=
\chi ^T  \sum _{j=0}^{\nu}  \| \Pi _j \|^2 \widetilde{Z}_{jj} \chi , \label{ttt5}
\end{align}

\vspace{-0.6cm}
and

\vspace{-1.0cm}
\begin{align}
 \sum _{j=0}^{\nu} \mathop{\textrm{He}} \left( \chi ^T \widetilde{N}_j w_j \right)&=
 \sum _{k=0}^{\nu} \left( \chi ^T N_k \sum _{j=0}^{\nu} c_{kj} \langle f,\Pi_j \rangle 
+ \sum _{j=0}^{\nu} \langle f,\Pi_j \rangle ^T c_{kj}N_k^T \chi \right)
 . \label{ttt6}
\end{align}
Equation (\ref{ttt4}) follows from (\ref{ttt5}) and (\ref{ttt6}).
 Now, it has to be shown only that $\widetilde{\Psi} \geq 0.$ By Schur complements, this is equivalent to
\begin{align*}
       {\clr \widetilde{\Phi} =}   \begin{bmatrix}
           \widetilde{Z}_{00}-\widetilde{N}_0 W^{-1} \widetilde{N}_{0}^{T} & \ldots & \widetilde{Z}_{0\nu}-\widetilde{N}_0 W^{-1} \widetilde{N}_{\nu}^{T} \\
           \vdots & \ddots & \vdots \\
           \widetilde{Z}_{\nu 0}-\widetilde{N}_{\nu} W^{-1} \widetilde{N}_{0}^{T} & \ldots & \widetilde{Z}_{\nu \nu}-\widetilde{N}_{\nu} W^{-1} \widetilde{N}_{\nu}^{T} \\
         \end{bmatrix}
        \geq 0.
\end{align*}
The non-diagonal blocks are zeros because of the definition of $\widetilde{Z}_{kl}$ for $k \neq l .$
The diagonal blocks can be written as
\begin{align*}
\widetilde{\phi}_{jj} &=\widetilde{Z}_{jj} - \widetilde{N}_j W^{-1} \widetilde{N}_j^T  
 = \sum _{k=0}^{\nu} \sum _{l=0}^{\nu} c_{kj}c_{lj} \left( Z_{kl}-N_k W^{-1}N_l^T \right) .
\end{align*}
Let {\clr
$\widetilde{y} \in \mathbf{R}^{M_{{\clr 2}}} $ be arbitrary and let ${y}_j= \mathop{\textrm{col}}\left\{c_{0j}\widetilde{y}, \ldots , \right.$ $ \left. c_{\nu j}\widetilde{y} \right\},$ then $\widetilde{y}^T \widetilde{\phi}_{jj} \widetilde{y} = {{y}_j}^T {\Phi} {y}_j $ holds true with matrix ${\Phi}$ defined by blocks ${\phi}_{kl}=Z_{kl}-N_k W^{-1}N_l^T $. By Schur complements, $\Psi \geq 0$ is equivalent to ${\Phi} \geq 0,$ therefore $\widetilde{\phi}_{jj} \geq 0$ for all $j,$
i.e. matrix $\widetilde{\Psi}$ is positive semi-definite} indeed.

Having regard to the equivalence of GFMB and IFB-GFMB, we will consider only GFMB in the sequel.

\subsection{Simplified GFMB and Bessel based inequalities}

In what follows we shall need the following notations.
For $i=1,2, $ set {\clr $\mathcal{W}_i = \mathop{\textrm{diag}} \left\{ \frac{1}{\rho_{0}^{i}},..., \frac{1}{\rho_{\nu}^{i}} \right\}\otimes W,$}
$\mathcal{W}_{i-} =  \mathcal{W}_{i}^{-1},$
  ${\widehat{\mathcal{W}}}_{1} = \mathop{\textrm{diag}} \left(  \mathcal{W}_{1} , 0  \right) ,$
$\widehat{\mathcal{W}}_{2} = \mathop{\textrm{diag}} \left( 0,  \mathcal{W}_{2}   \right), $ and $\widehat{\mathcal{W}}_{i-}$ is defined analogously with $\mathcal{W}_{i-}$.

\begin{cor}[S-GFMB]\label{cor:1}
If $W \in \mathbf{S}_{n}^{+} ,$
\begin{align}
 & \langle f_i , Wf_i \rangle  \geq
- \mathop{\emph{\textrm{He}}} \left(  \chi _{i}^{T} N^i \underline{w}  \right)
-   \chi _{i}^{T} N^i \widehat{\mathcal{W}}_{i-}  {N^{i}}^T \chi _{i} ,      \label{ttt10}
\end{align}
where $\chi ^{i}$ and $ N^i$ are the same as in Lemma \ref{lem:1}. Moreover, the right hand side of (\ref{ttt10}) is always greater or equal than the right hand side of (\ref{ttt1}).
\end{cor}
{\clr \textbf{Proof.} The right hand side of (\ref{ttt1}) can be written as
\begin{align}
- \chi _{i}^T \left( \sum _{k=0}^{\nu} \rho_{k}^{i} Z_{kk}^{i} \right) \chi _{i} &-  \mathop{{\textrm{He}}} \left(   \chi _{i}^{T} N^i \underline{w}, \right)
=   -  \chi _{i}^T \left( \sum _{k=0}^{\nu} \rho_{k}^{i}
\left( Z_{kk}^{i}-N_k^i W^{-1}   {N_{k}^{i}}^T  \right) \right) \chi _{i} \nonumber \\
&\hspace{0.60cm} - \mathop{\textrm{He}} \left(  \chi _{i}^T N^i \underline{w}  \right) 
-  \chi _{i}^T \left( \sum _{k=0}^{\nu}  N^i \widehat{\mathcal{W}}_{i-}  {N^{i}}^T  \right) \chi _{i}. \label{ttt9}
\end{align}
It follows from $\Psi^i \geq 0$  that
the choice of matrices
$Z_{kl}^{i}=N_k^i W^{-1}  {N_{0}^{i}}^T$ is admissible, thus (\ref{ttt9}) implies (\ref{ttt10}).
The second part follows from
$Z_{kk}^{i} \geq  N_k^i W^{-1}   {N_{k}^{i}}^T.$
 }
\begin{rem}\label{rem:3}
Corollary \ref{cor:1} is a straightforward extension of Lemma 5 of \cite{zha17} (continuous time case), and Lemma 2/6 of \cite{zha16b} (discrete time case).
\end{rem}
\begin{cor}[BBI]\label{cor:2}
If $W \in \mathbf{S}_{n}^{+} ,$
\begin{align}
 & \langle f_i , Wf_i \rangle  \geq {\underline{w}^{i}}^T \mathcal{W}_i \underline{w}^i.
    \label{c27_2}
\end{align}
Moreover, the right hand side of (\ref{c27_2}) is always greater or equal than the right hand side of (\ref{ttt1}) and (\ref{ttt10}).
\end{cor}
{\clr \textbf{Proof.}  The right hand side of (\ref{ttt10}) can be written as

\vspace{-1.1cm}
\begin{align}
&- \mathop{\textrm{He}} \left(   \chi _{i}^{T} N^i \underline{w}  \right)
-    \chi _{i}^{T} N^i \widehat{\mathcal{W}}_{i-}  {N^{i}}^T \chi _{i} \nonumber \\
= &-\left(  \chi _{i}^T N^i+\underline{w}^T \widehat{\mathcal{W}}_i
\right)
  \widehat{\mathcal{W}}_{i-}
         \left(  {N^{i}}^T \chi _i + \widehat{\mathcal{W}}_i
                                                  \underline{w} \right) 
                                                    +\underline{w}^T \widehat{\mathcal{W}}_i
             \underline{w} ={\underline{w}^{i}}^T \mathcal{W}^i \underline{w}^i , \label{c27_1}
\end{align}

\vspace{-0.4cm}
for $i=1,2.$
Since the first term in the second line of (\ref{c27_1}) is non-positive,
and it is equal to zero, if ${N^{i}}^T \chi _i= - \widehat{\mathcal{W}}_i  \underline{w}$
the statements of the corollary follow.
}
\begin{rem}\label{rem:4}
The estimation of Corollary \ref{cor:2} is the same as that of Lemma 1 of \cite{gyet16}, which implies among others the Bessel-Legendre inequality of \cite{seu14}, \cite{gyet16}, and the Bessel-Chebyshev inequality of \cite{gyet17}.
\end{rem}
The preceding results can be summarized as follows.
\begin{thm}\label{thm:2}
The estimations (\ref{c26_1})-(\ref{ttt1}) of Lemma \ref{lem:1}, (\ref{ttt10}) of Corollary \ref{cor:1} and (\ref{c27_2}) of Corollary \ref{cor:2}
are equivalent, and the tightest estimation is obtained under the choice of free parameters yielding (\ref{c27_2}).
\end{thm}

\subsection{Simplified FMB}

Choosing $\chi _i=\underline{w}$ we obtain the following estimation from Corollary \ref{cor:1}.
\begin{cor}[S-FMB]\label{cor:4}
If $W \in \mathbf{S}_{n}^{+} ,$

\vspace{-0.8cm}
\begin{align}
 & \langle f_i , Wf_i \rangle  \geq
- \underline{w}^T \left(\mathop{\emph{\textrm{He}}} \left(   N^i   \right)
-    N^i \widehat{\mathcal{W}}_{i-}  {N^{i}}^T \right) \underline{w} ,      \label{ttt100}
\end{align}

\vspace{-0.6cm}
where  $ N^i$ is the same as in Lemma \ref{lem:1}.
\end{cor}

\begin{rem}\label{rem:5}
Estimation (\ref{ttt100}) of Corollary \ref{cor:4}
is a general formulation of several FMB results. By appropriate choice of the matrices, of $\nu$
and of the polynomials, one obtains from S-FMB among others Lemma 3 of \cite{chen16b}, Theorem 1 of \cite{chen16c},
Lemma 4 of \cite{zeng15a}, Lemma 1 of \cite{zeng15b},
Lemma 2 of \cite{liu17}, Lemma 3 of \cite{seu14}, Lemma 4 of \cite{zha16a} (continuous time case)
and Lemma 2 of \cite{chen16a}, Corollaries 1-3 of \cite{chen16d}, Lemmas 7-9 of \cite{lee17b}, Lemmas 2-3 of \cite{wan16},
Lemma 2 of \cite{chen16a} (discrete time case).
\end{rem}
\begin{thm} \label{thm:4}
The estimations (\ref{ttt10}) of Corollary \ref{cor:1} and (\ref{ttt100}) of Corollary \ref{cor:4} are equivalent.
\end{thm}
\textbf{Proof.} It is enough to show that for any given $\chi _i$ and $N^i$ there is a $\widetilde{N}^i$ such that the right hand side of (\ref{ttt100}) equals to the right hand side of (\ref{ttt10}) taken with $N^i=\widetilde{N}^i.$ In fact, there exist an orthogonal matrix $Q$ and a scalar $\eta$ such that
$ \chi _i= \eta Q \underline{w}. $ Substituting it into (\ref{ttt10}) and taking $\widetilde{N}^i=\eta Q N^i$ results in (\ref{ttt100}).
{\clr
\begin{rem} \label{rem:6}
Choosing $\chi _i=\underline{w},$ one can derive from GFMB the FMB inequality, and its equivalence can be proven in the same way as Theorem \ref{thm:4}.
\end{rem}
Having regard the previous equivalences, we will consider only S-FMB and BBI in the sequel.
}

\section{Estimations for two connected intervals}

Specify now the Euclidean spaces  $\mathbb{E}_i$ $(i=0,1,2)$ according to Remark \ref{rem:0}. Denoting the lengths of the intervals as
$h=b-a,$ $h_1 =c-a,$ $h_2 =b-c$, where $a<c<b,$ one can verify that
\begin{align*}
& \rho_j^1 =  \frac{h_1}{2j+1} , \; \rho_j^2 = \frac{h_2}{2j+1},  \;
\rho_j^0 =  \frac{h}{2j+1}  ,
\end{align*}
satisfy condition $\rho _{k}^{i} \geq \| \Pi_{ki} \|_{i}^{2}.$
Introducing the notations 
$\alpha =\frac{h_1}{h} , $  {\clr $\beta =\frac{h_2}{h} $ $= 1-\alpha,$ 
$\mathcal{{W}}= \mathop{\textrm{diag}} \left\{ 1,3,...,2\nu +1 \right\} \otimes W , \;$
we} obtain from Corollary \ref{cor:4} that
\begin{align}
  \mathop{\textrm{(DS-FMB)}} \hspace{0.3cm}
\langle f , W f \rangle_{{\clr {0}}}
= \langle f_1 , W f_1 \rangle _1
+& \langle f_2 , W f_2 \rangle _2   \nonumber \\
  \geq & \frac{1}{h}\underline{w}^T   \Omega_F (\alpha ,\widehat{N}^1, \widehat{N}^2) \underline{w} , \label{ttt13}
\end{align}

\vspace{-0.7cm}
where

\vspace{-0.9cm}
\begin{align*}
\Omega_F (\alpha ,\widehat{N}^1, \widehat{N}^2)
& =
  \mathop{\textrm{He}} \left(- h \begin{bmatrix} \widehat{N}^1 & 0 \end{bmatrix}+\begin{bmatrix} 0 & \widehat{N}^2 \end{bmatrix}
  \right) \\
 & 
 + (-h
           \widehat{N}^1 )
                      (\alpha \mathcal{{W}}^{-1})
           (-h            \widehat{N}^1)^T
+  (-h
           \widehat{N}^2 )
                      (\beta \mathcal{{W}}^{-1})
           (-h            \widehat{N}^2)^T
                     .
\end{align*}
On the other hand, it follows from Corollary \ref{cor:2} that
\begin{align}
  \mathop{\textrm{(DBBI)}} \hspace{0.3cm}
\langle f , W f \rangle
&= \langle f_1 , W f_1 \rangle _1
+ \langle f_2 , W f_2 \rangle _2   
 \geq \frac{1}{h}\underline{w}^T   \Omega_B (\alpha ) \underline{w} ,  \label{c27_3}
\end{align}

\vspace{-0.9cm}
where

\vspace{-0.9cm}
\begin{align*}
& \Omega_B (\alpha) =
                          \begin{bmatrix}
                            \frac{1}{\alpha}\mathcal{{W}} & 0 \\
                            0 & \frac{1}{\beta}\mathcal{{W}} \\
                          \end{bmatrix}
                               .
\end{align*}
From the results of the previous subsections, it follows that (\ref{c27_3}) yields a tighter lower bound then (\ref{ttt13}). However,
if (\ref{c27_3}) is applied to the stability analysis of systems with time-varying delays, the obtained estimation is non-convex in the lengths of the intervals. To avoid the non-convexity, \emph{one has to apply a further lower estimation of the right hand side of (\ref{c27_3})}. Therefore, a real comparison can be obtained by relating (\ref{ttt13}) to the combination of (\ref{c27_3}) and some convexifying lower bound.

In what follows, we shall investigate the application of the classical reciprocally convex combination (RCC) lemma of \cite{park11}, the extended reciprocally convex combination (ERC) lemma of \cite{seu16}, and Lemma 3 of \cite{liuk16}, which 
will be referred to as (M-LSX). We shall apply the (ERC) lemma with two special choices of the matrices as well.
These estimations are summarized in the lemma below.

\begin{lem} \label{lem:10}
If $W \in \mathbf{S}_{n}^{+} ,$ then for $k=1, \ldots,5,$
\begin{align}
\Omega _B (\alpha)
 \geq  \Omega _k (\alpha), \; &\forall \alpha \in (0,1), \label{c0301_1} 
\nonumber
\end{align}
where $\Omega _k$ is given in different cases 
as follows:
\begin{itemize}
\item  \emph{(M-LSR):} with arbitrary $V_1, V_2 \in \mathbf{R}^{{\clr M_2} \times M_1},$
\begin{align}
 \Omega_1(\alpha) = \Omega_1(\alpha, V_1,V_2) & = \mathop{\textrm{\emph{{He}}}}\left(
  \begin{bmatrix}
                                    V_1 & 0 \\
                                  \end{bmatrix}
                                  +   \begin{bmatrix}
                                    0 & V_2 \\
                                  \end{bmatrix}   \right) \nonumber\\
  & -\alpha V_1 \mathcal{{W}}^{-1}V_1^T-\beta V_2 \mathcal{{W}}^{-1} V_2^T;  
  \end{align}
\item \emph{(ERC):} with arbitrary $  Y_1, Y_2 \in \mathbf{R}^{{\clr M_2} \times M_1}$ and $X_1, X_2 \in \mathbb{S}_{M_1}^{+},$  satisfying
 for $ \alpha = 0, 1,$
\begin{align}
& \begin{bmatrix}\mathcal{{W}}& 0 \\ 0 & \mathcal{{W}} \end{bmatrix}- \alpha  \begin{bmatrix}X_1 & Y_1 \\ Y_1^T & 0 \end{bmatrix}
 -\beta \begin{bmatrix}0 & Y_2 \\ Y_2^T & X_2 \end{bmatrix} \geq 0, \label{SE1} 
\end{align}

\vspace{-0.5cm}
\begin{align}
   \Omega_2(\alpha) & = \Omega_2(\alpha,X_1,X_2,Y_1,Y_2) 
   = \begin{bmatrix}
                                    \mathcal{{W}}+\beta X_1 & \alpha Y_1+\beta Y_2 \\
                                    \ast & \mathcal{{W}}+\alpha X_2 \\
                                  \end{bmatrix} ; \label{c03013}
\end{align}
\item \emph{(SERC):} with arbitrary $  Y_1, Y_2 \in \mathbf{R}^{M_1 \times M_1}$ and  $\widehat{X}_1= \mathcal{{W}}-Y_1 \mathcal{{W}}^{-1} Y_1^T,$ $\widehat{X}_2= \mathcal{{W}}-Y_2^T \mathcal{{W}}^{-1} Y_2,$
    \begin{align}
   \Omega_3(\alpha) & = \Omega_3(\alpha,Y_1,Y_2) 
    =\begin{bmatrix}
                                    \mathcal{{W}}+\beta \widehat{X}_1 & \alpha Y_1+\beta Y_2 \\
                                    \ast & \mathcal{{W}}+\alpha \widehat{X}_2 \\
                                  \end{bmatrix}; \label{c03014}
\end{align}
\item \emph{(MERC):} with arbitrary $  Y \in \mathbf{R}^{M_1 \times M_1}$ and  $\overline{X}_1= \mathcal{{W}}-Y \mathcal{{W}}^{-1} Y^T,$ $\overline{X}_2= \mathcal{{W}}-Y^T \mathcal{{W}}^{-1} Y,$
\begin{align}\label{c03015}
   \Omega_4(\alpha) & =  \Omega_4(\alpha,Y) =
   \begin{bmatrix}
                                    \mathcal{{W}}+\beta\overline{X}_1 &  Y \\
                                    \ast & \mathcal{{W}}+\alpha \overline{X}_2 \\
                                  \end{bmatrix};
\end{align}
\item \emph{(RCC):} with arbitrary $  Y \in \mathbf{R}^{M_1 \times M_1}$ satisfying (\ref{SE1}) with $X_1=X_2=0$ and $Y_1=Y_2=Y,$
\begin{align}
 &\Omega_5(\alpha)=   \Omega_5(\alpha,Y) = \begin{bmatrix}\mathcal{{W}}& Y \\ Y & \mathcal{{W}} \end{bmatrix} . \label{c03016}
\end{align}
\end{itemize}
\end{lem}
\begin{thm}\label{thm:3}
The estimation (\ref{ttt13}) and the
estimations of Lemma \ref{lem:10} are related as follows:
\begin{itemize}
 \item
[\emph{\bf{(A)}}] \emph{(DS-FMB)}
is equivalent to  \emph{(DBBI)} $\&$ \emph{(M-LSR)};
 \item
[\emph{\bf{(B)}}] \emph{(DBBI)} $\&$  \emph{(M-LSR)} implies \emph{(DBBI)} $\&$  \emph{(SERC)},
but not conversely;
 \item
[\emph{\bf{(C)}}] \emph{(DBBI)} $\&$  \emph{(ERC)}
is equivalent to  \emph{(DBBI)} $\&$  \emph{(SERC)};
 \item
[\emph{\bf{(D)}}] \emph{(DBBI)} $\&$  \emph{(SERC)} implies \emph{(DBBI)} $\&$  \emph{(MERC)},
but not conversely;
 \item
[\emph{\bf{(E)}}] \emph{(DBBI)} $\&$  \emph{(MERC)} implies \emph{(DBBI)} $\&$  \emph{(RCC)},
and $\Omega _5(\alpha,Y) \leq \Omega _4(\alpha,Y) $ if $Y$ is chosen according to (RCC).
\end{itemize}
\end{thm}
\textbf{Proof.} The proof of \textbf{(A)} consists of verifying that  
$$\Omega _B(\alpha,\widehat{N}_1,\widehat{N}_2)=\Omega _1(\alpha,-h \widehat{N}_1, -h \widehat{N}_2).$$

In order to prove \textbf{(B)}, we can verify first that the choice of
\begin{equation*}
  V_1=\begin{bmatrix}
                                    \mathcal{{W}} \\
                                    Y_2^T \\
                                  \end{bmatrix}, \; \mbox{ and }
 V_2=\begin{bmatrix}
                                   Y_1  \\
                                  \mathcal{{W}}   \\
                                  \end{bmatrix},
\end{equation*}
{\clr yields} $\Omega _1(\alpha)=\Omega _3(\alpha),$ i.e. $\Omega _3(\alpha)$ can be obtained as a special case of $\Omega _1(\alpha).$ To show the
second part of the statement, let us write $V_1=\mathop{\textrm{col}} \left\{ \mathcal{{W}}+\Xi _1, V_{12} \right]$ and
$V_2=\mathop{\textrm{col}} \left\{ V_{21}, \mathcal{{W}}+\Xi _2 \right].  $  
Then

\vspace{-0.6cm}
\begin{align*}
  \Omega _1(\alpha,V_1,V_2) & -\Omega _3(\alpha,Y_1,Y_2)  \\
 & = \begin{bmatrix}  -\Xi_1 \mathcal{{W}}^{-1} \Xi_1^T  & \Xi_1 \mathcal{{W}}^{-1} V_{12}^T + V_{12}-Y_1  \\
 \ast & \mathop{\textrm{He}}(\Xi_2)  -V_{12} \mathcal{{W}}^{-1} V_{12}^T +Y_1^T \mathcal{{W}}^{-1}Y_1\end{bmatrix}.
\end{align*}
If $\Xi_1\neq 0,$ then there is a $y_1 \in \mathbf{R}^{M_1}$ such that $y_1^T\left( \Omega _1(\alpha)_{11}-\Omega _3(\alpha)_{11} \right) y_1 < 0 $ for $\alpha=1,$ and for continuity reasons, for  $\alpha \in (1-\delta, 1], (0<\delta<1)$
\emph{independently of the choice of }$Y_1, \ Y_2.$
At the same time, if $\Xi_2$ is such that $\Xi_2+\Xi_2^T   -V_{12} \mathcal{{W}}^{-1} V_{12}^T>0,$
then there is a $y_2 \in \mathbf{R}^{M_1}$ such that $y_2^T\left( \Omega _1(1)_{22}-\Omega _3(1)_{22} \right) y_2 > 0 $
independently of the choice of $Y_1, \ Y_2.$
Therefore,
$
 \Omega _1(\alpha,V_1,V_2)-\Omega _3(\alpha,Y_1,Y_2)     
$
may be indefinite for some $V_1, V_2$ independently of the choice of $Y_1, \ Y_2.$ This means that neither $\Omega _1(\alpha)\geq \Omega _3(\alpha),$ nor $\Omega _1(\alpha)\leq \Omega _3(\alpha)$ is true in general.

To show \textbf{(C)} we observe that it follows from (\ref{SE1}) that ${X}_1 \leq \mathcal{{W}}-Y_1 \mathcal{{W}}^{-1} Y_1^T,$ and ${X}_2= \mathcal{{W}}-Y_2^T \mathcal{{W}}^{-1} Y_2.$ Thus $\Omega _3(\alpha)$ is a special case of $\Omega _2(\alpha),$ and
$\Omega _2(\alpha) \leq \Omega _3(\alpha)$ for any $Y_1, Y_2,$ which proves \textbf{(C)}.

{\clr Since $\Omega _4(\alpha,Y)=\Omega _3(\alpha,Y,Y),$ the first assertion \textbf{(D)} is obvious.}  To show the
second part of the statement, let us write $Y_i=Y+\Upsilon _i \ (i=1,2).$ Then a straightforward computation shows that $\Omega _3(1,Y+\Upsilon _1,Y+\Upsilon _2)-\Omega _4(1,Y)$ is indefinite, if $\Upsilon _1 \neq 0, $ thus neither $\Omega _3(\alpha)\geq \Omega _4(\alpha),$ nor $\Omega _3(\alpha)\leq \Omega _4(\alpha)$ is true.

The proof of \textbf{(E)} is immediate, since $\Omega _5(\alpha,Y)=\Omega _4(\alpha,Y){\clr ,}$ 
if $Y\in \mathbf{R}^{M_1 \times M_1}$ is such that (\ref{SE1}) is satisfied with $X_1=X_2=0,$  $Y_1=Y_2=Y,$ and {\clr $\overline{X}_1=\overline{X}_2=0$} is substituted in $\Omega _4$.
Since $\overline{X}_1, \  \overline{X}_2\geq 0$ for the above
$Y{\clr ,} $ inequality
$\Omega _5(\alpha) \leq \Omega _4(\alpha) $ is true.
\begin{rem}\label{rem:7}
1.) Lemma 4 of \cite{zha16a} is a special case of the estimation (DBBI) $\&$ (MERC).

2.)
 We have seen that neither of the estimations obtained by (M-LSR) and (SERC) is better then the other. Nevertheless, (M-LSR) may be advantageous in the analysis of systems with time-varying delays. Similar note is due with respect to the estimations (SERC) and (MERC) on favour of (SERC).
\end{rem} 

\section{Conclusions}

We introduced an extension of the generalized free matrix based inequality in a unified form suitable for the estimation of integrals and sums of quadratic functions. The equivalences of several known variants were proven, including the {\clb free-matrix-based} inequalities and its simplified form. It was shown that the Bessel-based estimation is is at least as good as any of the others, while the S-FMB estimation is at least as good as GFMB, IFB-GFMB. Secondly, the relationship between the S-FMB estimation and the combination of the Bessel-based inequality with different bounding inequalities being affine in the length of the intervals were intensivly investigated.



\begin{thebibliography}{99}     
\bibitem[Chen et~al. (2016a)]{chen16a} Chen J., Lu J. \& Xu S. (2016a). Summation inequality and its application
to stability analysis for time-delay systems. \emph{IET Control Theory \& Applications, 10,} 391-395.

\bibitem[Chen et~al. (2016b)]{chen16b} Chen J., Xu S., Zhang B. \& Liu G. (2016b). A note on relationship between two classes of
integral inequalities. \emph{IEEE Transactions on Automatic Control}, doi: 10.1109/TAC.2016.2618367.

\bibitem[Chen et~al. (2016c)]{chen16c} Chen J., Xu S. \& Zhang B. (2016c). Single/multiple integral inequalities with
applications to stability analysis of time-delay systems. \emph{IEEE Transactions on Automatic Control};, doi: 10.1109/TAC.2016.2617739.

\bibitem[Chen et~al. (2016d)]{chen16d} Chen J., Xu S., Jia X. \& Zhang B. (2016d). Novel summation inequalities and their
applications to stability analysis for systems with time-varying delay.
\emph{IEEE Transactions on Automatic Control}, doi: 10.1109/TAC.2016.2606902.

\bibitem[Chen et~al. (2016e)]{chen16e}  Chen J., Xu S., Chen W., Zhang B.,  Ma Q. \&
 Zou Y. (2016e). Two general integral inequalities and their applications to stability
analysis for systems with time-varying delay. \emph{International Journal of Robust Nonlinear Control, 26,} 4088-4103.

\bibitem[Gyurkovics (2015)]{gye15} Gyurkovics \'{E}. (2015). A note on Wirtinger-type inequalities for time-delay systems.
\emph{Automatica, 61,} 44-46.

\bibitem[Gyurkovics et~al. (2016)]{gyet16} Gyurkovics \'{E.} \& Tak\'{a}cs T. (2016). Multiple integral inequalities and stability analysis
of time delay systems. \emph{Systems \& Control Letters, 96,} 72-80.

\bibitem[Gyurkovics et~al. (2017)]{gyet17} Gyurkovics \'{E}., Kiss K., Nagy I. \& Tak\'{a}cs T. (2017). Multiple summation inequalities and
their application to stability analysis of discrete-time delay systems. \emph{Journal of the Franklin Institute, 354,} 123-144.


\bibitem[Hien et~al. (2015)]{hien15} Hien L.V. \& Trinh H. (2015). An enhanced stability criterion for time-delay systems
via a new bounding technique. \emph{Journal of the Franklin Institute, 352,} 4407-4422.

\bibitem[Lee et~al. (2015)]{leejfi15}  Lee T.H.,  Park J.G.,  Park M.J.,  Kwon O.M. \&  Jung H.-Y. (2015).
On stability criteria for neural networks with time-varying delay using
Wirtinger-based multiple integral inequality,
\emph{Journal of the Franklin Institute, 352,} 5627-5645.

\bibitem[Lee et~al. (2017a)]{lee17a} Lee S.Y., Lee W.I. \& Park P. (2017a). Polynomials-based integral inequality for stability analysis of linear systems with time varying delays.
    \emph{Journal of the Franklin Institute, 354,} 2053-2067.

\bibitem[Lee et~al. (2017b)]{lee17b} Lee S.Y., Lee W.I. \& Park P. (2017b). Polynomials-based summation inequalities and their applications to discrete-time systems with time-varying delays.
    \emph{International Journal of Robust and Nonlinear Control}, doi: 10.1002/rnc.3755.



\bibitem[Liu et~al.(2016)]{liuk16} Liu K., Seuret A. \& Xia Y. (2016). Stability analysis of systems with time-varying delays
via the second-order Bessel-Legendre inequality. \emph{Automatica, 76,} 138-142.

\bibitem[Liu et~al.(2017)]{liu17} Liu Y., Park J.H. \& Guo B.-Z. (2017). New results on stability of linear systems with time varying
delay. \emph{IET Control Theory \& Applications, 11,} 129-134.


\bibitem[Park et~al.(2011)]{park11}  Park P.G.,   Ko J.W. \&  Jeong C. (2011).   Reciprocally convex approach to stability of
systems with time-varying delays. \emph{Automatica, 47,}(1) 235-238.

\bibitem[Seuret et~al.(2013)]{S-GAut13}  Seuret A. \&  Gouaisbaut F. (2013).  Wirtinger-based integral inequality: Application to time-delay systems. \emph{Automatica, 49,}  2860-2866.

\bibitem[Seuret et~al.(2014)]{seu14} Seuret A. \& Gouaisbaut F. (2014). Complete quadratic Lyapunov functionals using
Bessel-Legendre inequality. \emph{European Control Conference};  Strasbourg, France, (pp. 448-453).

\bibitem[Seuret et~al.(2016)]{seu16} Seuret A. \& Gouaisbaut F. (2016). \emph{Delay-dependent reciprocally convex combination lemma}.
Rapport LAAS n16006. 2016. $\langle$hal-01257670$\rangle$

\bibitem[Wan et~al.(2016)]{wan16} Wan X., Wu M., He Y. \& She J. (2016). Stability analysis for discrete time-delay systems based on new
finite-sum inequalities. \emph{Information Sciences, 369,} 119-127.

\bibitem[Zeng et~al.(2015a)]{zeng15a} Zeng H.B., He Y., Wu M. \& She J. (2015a). Free-matrix-based integral inequality for stability
analysis of systems with time-varying delay. \emph{IEEE Transactions on Automatic Control, 60,} 2768-2772.

\bibitem[Zeng et~al. (2015b)]{zeng15b} Zeng H.B., He Y., Wu M. \& She J.H. (2015b). New results on stability analysis
for systems with discrete distributed delay. \emph{Automatica, 60,} 189-192.

\bibitem[Zhang et~al.(2016a)]{zha16a} Zhang C.K., He Y., Jiang L., Wu M. \& Zeng H.B. (2016a). Stability analysis of systems with time-varying
delay via relaxed integral inequalities. \emph{Systems \& Control Letters, 92,} 52-61.

\bibitem[Zhang et~al.(2016b)]{zha16b} Zhang C.K., He Y., Jiang L., Wu M. \& Zeng H.B. (2016b). Summation inequalities to bounded real lemmas of discrete-time systems with time-varying delay. \emph{IEEE Transactions on Automatic Control}, doi: 10.1109/TAC.2016.2600024.

\bibitem[Zhang et~al.(2017)]{zha17} Zhang C.K., He Y., Jiang L., Lin W.J. \& Wu M. (2017).
Delay-dependent stability analysis of neural networks with time-varying delay: A generalized free-weighting-matrix approach.
\emph{Applied Mathematics and Computation, 294,} 102-120.
\end{thebibliography}


\end{document}